\documentclass[12pt,a4paper]{article}

\usepackage{latexsym}

\newcommand{\be}{\begin{equation}}	
\newcommand{\ee}{\end{equation}}	
\newcommand{\bern}{\begin{eqnarray*}}	
\newcommand{\eern}{\end{eqnarray*}}	
\newcommand{\beqp}{\begin{eqproof}}	
\newcommand{\eeqp}{\end{eqproof}}	

\newcommand{\bt}{\begin{teorema}}	
\newcommand{\et}{\end{teorema}} 	
\newcommand{\bd}{\begin{definizione}}	
\newcommand{\ed}{\end{definizione}}	
\newcommand{\bc}{\begin{corollario}}	
\newcommand{\ec}{\end{corollario}}	
\newcommand{\bp}{\begin{dimostrazione}}	
\newcommand{\ep}{\end{dimostrazione}}	
\newcommand{\bl}{\begin{lemma}}		
\newcommand{\el}{\end{lemma}}		
\newcommand{\bpr}{\begin{proposizione}}	
\newcommand{\epr}{\end{proposizione}}	
\newcommand{\besi}[1][{}]
	{\begin{esempi}\rm\textbf{#1}\begin{enumerate}}	
\newcommand{\eesi}{\end{enumerate}\end{esempi}}	
\newcommand{\beso[1]}{\begin{esempio}[#1]\rm}	
\newcommand{\eeso}{\end{esempio}}	
\newcommand{\boss}{\begin{osservazione}\rm}	
\newcommand{\eoss}{\end{osservazione}}	
\newcommand{\refe}[1]{(\ref{#1})}	

\font\scdoppio=msbm8

\newcommand{\C}{\mbox{\corsivo C}}	
\newcommand{\F}{\mbox{\corsivo F}}	
\newcommand{\K}{\mbox{\corsivo K}}	
\newcommand{\R}{\hbox{\doppio R}}	
\newcommand{\RR}{\hbox{\scdoppio R}}	
\newcommand{\nat}{\hbox{\doppio N}}	
\newcommand{\uni}{\mathbf{1}}		

\newcommand{\norm}[1]{\left\|#1\right\|}		
\newcommand{\abs}[1]{\left|#1\right|}	
\newcommand{\decl}{:=}			

\renewcommand{\theequation}{\thesection.\arabic{equation}}

\newtheorem{teorema}{Theorem}[section]
\newtheorem{definizione}[teorema]{Definition}
\newtheorem{corollario}[teorema]{Corollary}
\newtheorem{lemma}[teorema]{Lemma}
\newtheorem{proposizione}[teorema]{Proposition}
\newtheorem{esempio}[teorema]{Example}
\newtheorem{esempi}[teorema]{Examples}
\newtheorem{osservazione}[teorema]{Remark}

\newcounter{proofeqno}		
\newenvironment{dimostrazione}[1]	
	{\noindent \textbf{Proof #1.\ }\setcounter{proofeqno}{0}}
	{\hfill $\Box$\par\medskip}

\newenvironment{eqproof}	
	{\addtocounter{proofeqno}{1}$$}
	{\eqno(\theproofeqno)$$}

	{\addtocounter{proofeqno}{1} \addtocounter{equation}{-1}
		\renewcommand\theequation{\theproofeqno}
		\begin{eqnarray}   }
	{ \end{eqnarray}
		\renewcommand\theequation{\thesection.\@arabic\c@equation}	}
	{\mbox{}\par }{\par\medskip}

\def\acknowledgment{\goodbreak\subsection*{Acknowledgment}
\bgroup \footnotesize}

\def\acknowledgments{\goodbreak\subsection*{Acknowledgments}
\bgroup \footnotesize}

\def\endacknowledgment{\vskip1sp\egroup}
\def\endacknowledgments{\vskip1sp\egroup}

\newcommand{\pr}{\mathrm{pr}}
\newcommand{\vi}{{\bf v}}

\font\corsivodo=rsfs12 
\font\doppiodo=msbm12 

\renewcommand{\C}{\mbox{\corsivodo C}}	
\renewcommand{\F}{\mbox{\corsivodo F}}	
\renewcommand{\K}{\mbox{\corsivodo K}}	
\renewcommand{\R}{\hbox{\doppiodo R}}	
\newcommand{\natur}{\hbox{\doppiodo N}}	
\renewcommand{\RR}{\hbox{\scdoppio R}}	
\renewcommand{\nat}{\hbox{\scdoppio N}}	

\begin{document}
\title{Extension of Bernstein Polynomials to Infinite Dimensional Case}
\author{Lorenzo D'Ambrosio \\
  {\small Dipartimento di Matematica, Universit\`a di Bari,}\\
  {\small  via Orabona, 4\ \   I--70125 Bari Italy,}\\ 
  {\small e-mail \tt dambros@dm.uniba.it}}

\date{}

\maketitle\label{ver 2.2}

{\abstract The purpose of this paper is to study some new concrete
approximation processes for continuous vector-valued mappings defined
on the infinite dimensional cube or on a subset
of a real Hilbert space. 
In both cases these operators are modelled on classical
Bernstein polynomials and represent a possible extension to an infinite
dimensional setting.

The same idea is generalized to obtain from a given approximation process
for function defined on a real interval
a new approximation process
for  vector-valued mappings defined on subsets of a real Hilbert space.



}

\section{Introduction} 

The purpose of this paper is to define an explicit sequence
of operators that is an approximation process for continuous vector-valued
mappings $F:X\to E$, where $X$ has ``infinite dimension''.
More precisely we deal with two cases.
The first is when $X$ is the cube
\[ C_\infty\decl [0,1]^{{\nat}^*},\]
with the canonical product topology, where $\natur^*$ denotes the set 
$\natur\setminus\{0\}$ and $\natur\decl\{0,1,2,\dots\}$.
The other case we consider, is when $X$ is an unbounded, closed subset of
a real Hilbert space endowed with the weak topology.

A first approximation process, we are going to construct, is modelled
on the Bernstein polynomials. Later we shall give a generalization of this construction.

The Bernstein polynomials, for a continuous function $F\in\C(C_k)$,  on the
$k$-dimensional cube  $C_k\decl[0,1]^k$, are defined
at $t=(t_1,\dots,t_k)\in C_k$, as
\[ B_{n,k}(F)(t)\decl \sum^n_{\stackrel{\stackrel{\scriptstyle j_1=0}{\dots}}{\scriptstyle j_k=0}}
	F\left(\frac{j_1} n,\dots, \frac {j_k}{n}\right) \psi_{n,j_1}(t_1)
	\cdots\psi_{n,j_k}(t_k), \label{bernstein}\]
where
\[\psi_{n,j}(t)\decl{n\choose j}t^{j}(1-t)^{n-j}.\]
It is well known that the sequence $(B_{n,k})_{n\ge 1}$ realizes an 
approximation process on $\C(C_k)$ as specified by
\bt\label{tb} \begin{enumerate} \item For any $F\in \C(C_k)$, 
	$B_{n,k}(F)\rightarrow F$ uniformly on $C_k$ as $n\rightarrow \infty$.
\item Let $C_k$ endow with the distance $d(x,y)\!\decl\!\sum_{i=1}^k\!\abs{x_i-y_i}$.
	If $F\!\in\! Lip_M(C_k)$, then
	$B_{n,k}(F)\in Lip_M(C_k)$.\footnote{Let $(X,d)$ be a metric space and
	$E$ normed space. A function $f:X\to E$ belongs to $Lip_M(X)$,
	if $\norm{f(t)-f(\tau)}\le Md(t,\tau)$, for any $t,\tau\in X$.}

\item For any convex function $F\in \C(C_k)$, $B_{n,k}(F)$ is convex
	with respect to each variable.

\item  For any $F\in \C(C_k)$, convex with respect to each variable and
	$n\ge1$, it results $F\le B_{n,k}(F)$.
\item  For any $F\in \C(C_k)$, convex with respect to each variable and
	$n\ge1$, it results $B_{n+1,k}(F)\le B_{n,k}(F)$.
\end{enumerate}
\et
We refer the interested reader to e.g. \cite{ac,dm,rasa94}.

Our  idea is simple. We link the index $n$ to the dimension $k$ of the
cube where the operator $B_{n,k}$ samples the function,
obtaining the operator $B_{n,n}$;
in the $C_\infty$ case,
the $n$-th operator acts sampling $F:C_\infty\to E$ on a $n$-dimensional cube.

In the next section we present the results, while the proofs are in
 the section~3.
The last section is devoted to extend the idea to other operators.

\section{Definitions and Results}
Let $X$ be a Hausdorff space and $E$ a normed space. We denote with $\F(X,E)$ and with
$\C(X,E)$ respectively the space
of all mappings $F:X\to E$ and its  subspace containing only the continuous mappings. 

Fix $g:X\to \R_+$, the symbol $\F(X,E,g)$ stands for the subspace of all
mappings $F$ belonging to $\F(X,E)$ such that $F/g$ is bounded.

For every $n\ge1$, we set
\[ A_n\decl\{h=(h_j)_{j\ge1} | h_j\in\natur,\ 0\le h_j\le n\ \mathrm{ for }\ j\le n,\ h_j=0\ \mathrm{ for }\ j>n \}. \]
In other words, $h\in A_n$ if and only if it has the form
$h=(h_1,\dots, h_n,0,0,\dots)$ with  $0\le h_j\le n$ for every natural
$j\in\{1,2,\dots ,n\}$.

\subsection{$C_\infty$ case}
As we have mentioned in the introduction, the topology in $C_\infty = [0,1]^{{\nat}^*}$
is the canonical product one; every point $t\in C_\infty$ is identified with the
sequence $(t_j)_{j\ge 1}$.

Let $n\ge1$ be natural number, $h=(h_j)_{j\ge1}\in A_n$.
Define the function
\be \varphi_{n,h}(t)\decl {n\choose h_1} \cdots{n\choose h_n}t_1^{h_1}(1-t_1)^{n-h_1}
	\cdots t_n^{h_n}(1-t_n)^{n-h_n},\label{weight}
\ee
for every $t\in C_\infty$. Notice that $\varphi$ has the form
$\varphi_{n,h}(t)=\psi_{n,h_1}(t_1)\cdots\psi_{n,h_n}(t_n)$.

For every $n\ge1$, $F:C_\infty\to E$ and $t\in C_\infty$, we define
\[ L^1_n(F)(t)\decl \sum_{h\in A_n} F\left(\frac h n\right)\varphi_{n,h}(t),\]
or, explicitly,
\bern L^1_n(F)(t) &=& \sum^n_{\stackrel{\stackrel{\scriptstyle h_1=0}{\dots}}{\scriptstyle h_n=0}}
	F\left(\frac{h_1} n,\dots, \frac{h_n} n,0,0,\dots\right)
	{n\choose h_1}t_1^{h_1}(1-t_1)^{n-h_1}\cdots\times\\
	&&\qquad\qquad\qquad\times\cdots{n\choose h_n}t_n^{h_n}(1-t_n)^{n-h_n}.
\eern

In section 3 we shall prove the following approximation result:
\bt\label{cube} For any $F\in\C(C_\infty,E)$, the convergence
	\[L^1_n(F)\rightarrow F\quad\mathrm{as}\quad n\rightarrow \infty \]
	holds uniformly on $C_\infty$.
\et

\subsection{Hilbert case}
Let $H$ be an infinite dimension separable real Hilbert space.
With $(a_j)_{j\ge1}$ we denote a Hilbert base of $H$, so that the points
$t\in H$ are represented by $t=\sum_{j=1}^\infty t_ja_j$.
A  well known fact says that $H$ is isometrically isomorphic to
the Hilbert space
$\ell^2\decl\{(t_n)_{n\ge1}|\sum_{n=1}^\infty |t_n|^2 < \infty\}$.
Therefore,  we shall use the  identification $H=\ell^2$.

We set
\[\Gamma \decl\{t\in H|\ 0\le t_i\le 1\}.\]

The definition of $\varphi_{n,k}$ in \refe{weight} is still valid for $t\in \Gamma$,
hence for every $n\ge1$, $F:\Gamma\to E$ and $t\in \Gamma$, we define
\be L^2_n(F)(t)\decl \sum_{h\in A_n} F\left(\frac h n\right)\varphi_{n,h}(t),\label{dl2} \ee
or, equivalently:
\bern L^2_n(F)(t) &=& \sum^n_{\stackrel{\stackrel{\scriptstyle h_1=0}{\dots}}{\scriptstyle h_n=0}}
	F\left(\frac{h_1} n,\dots, \frac{h_n} n,0,0,\dots\right)
	{n\choose h_1}t_1^{h_1}(1-t_1)^{n-h_1}\cdots\times\\
	&&\qquad\qquad\qquad\times\cdots{n\choose h_n}t_n^{h_n}(1-t_n)^{n-h_n}.
\eern

We remind the following definitions
\bd Let $X$ be a convex subset of a Banach space $Y$.
\begin{enumerate}\item The symbol $UCB(X,E)$ stands for the subspace
	of $\F(X,E)$ of all the uniformly continuous and bounded mappings.
	For $F\in UCB(X,E)$, we define, as usual, its
	\emph{modulus of continuity}, as
	\[\omega(F,\delta)\decl \sup\{\norm{F(u)-F(t)} | u,t\in X, \norm{u-t}\le \delta\} \qquad (\delta >0).\]

\item We say that $F\colon X\to E$ is
	\emph{weak-to-norm continuous} if it is continuous
	from $X$ equipped with the weak topology $\sigma(Y,Y')$ in $Y$,
	into $E$ with the norm topology.
	By $\K(X,E)$ we denote the space of all weak-to-norm continuous
	mappings from $X$ into $E$.
	We set $\K(X,E,g)\decl \K(X,E)\cap \F(X,E,g)$.
\end{enumerate}
\ed

The approximation results in the Hilbert case are as follows.

\bt\label{hilbert} 
	For any $F\in\K(\Gamma,E,1+\norm{\cdot}^2)$, the convergence
	\[L^2_n(F)(t)\rightarrow F(t)\]
	holds for any $t\in\Gamma$ and uniformly on relatively
	compact subsets of $\Gamma$.
\et

\bt\label{ucb} For any $F\in UCB(\Gamma,E)$, we have
	$L_n^2(F)\rightarrow F$ (as $n\rightarrow \infty$), uniformly
	on relatively compact subsets of $\Gamma$.
	Moreover for any $t\in \Gamma$, $n\ge 1$
	and $\delta >0$, there holds the estimate
	\[ \norm{L_n^2(F)(t)-F(t)}\le\omega(F,\delta)\left[1+\delta^{-2}
	(\sum_{j>n}t_j^2+\sum_{j=1}^n\frac{t_j-t_j^2}{n})\right],\]
	therefore, in particular
\bern \norm{L_n^2(F)(t)-F(t)}&\le&2\omega(F,\sqrt{\sum_{j>n}t_j^2+
		\sum_{j=1}^n\frac{t_j-t_j^2}{n}})\\
	&\le&2\omega(F,\sqrt{\sum_{j>n}t_j^2+\frac{\norm t}{\sqrt n}
		+\frac{\norm t^2}{n}}).\eern
\et

These operators $L_n^2$ satisfy the following preserving properties.
\bpr \label{lip}\begin{enumerate}
\item If $F\in Lip_M(\Gamma)$, then 
	$L^2_n(F)\in Lip_{\sqrt n M}(\Gamma)$ for any $n\ge 1$.
\item If $F\in\C(\Gamma,\R)$ is convex, then for any $n\ge1$,
	$L_n^2(F)$ is convex with respect to each variable.
\end{enumerate}
\epr

Thus, the analogues of the properties 1., 2. and 3. of Theorem \ref{tb}
are in some sense inherited from $L_n^2$.

Let $E$ be an ordered space. 
The following question arise.
What happens to properties 4. and 5.? They fail even in the case $E=\R$.
We shall prove this claim in the next section finding a counterexample.

\section{Proofs}

Before proving the statements of the previous section, we  recall
the following  definitions (cfr. \cite{prolla}).
For any function $g\in\F(X,\R)$ and any vector $\vi\in E$, with
$g\otimes \vi$ we denote the function belonging to $\F(X,E)$ defined as
$$(g\otimes \vi)(t)\decl g(t)\vi\qquad\mathrm{ for\ any}\quad t\in X.$$

\bd Let $S$ be a linear operator on $\F(X,\R)$. A linear operator $L$ on
	$\F(X,E)$ is said to be $S$-{\em regular} if
	\[L(g\otimes \vi)=S(g)\otimes \vi, \mathrm\ {for\ all\ }g\in\F(X,\R)\mathrm{\ and\ }\vi\in E.\]
	$L$ is said {\em monotonically regular}, if it is $S$-regular for some
	positive linear  operator on $\F(X,\R)$.
\ed

\boss\label{inter} The operators $L_n^1$ and $L_n^2$ are well
	defined on  scalar functions as well as on vector-valued mappings
	and we shall use the same symbol for the operators acting on
	vector-valued mappings or on scalar functions.
	Moreover, it is easily seen that both operators are monotonically regular.
\eoss

\subsection{Proof of Theorem  \ref{cube}}
Combining the results \cite[Theorem 4.4.6]{ac} and \cite[Theorem 9, pag. 111]{prolla}
we obtain
\bt\label{apt} Let $X$ be a compact Hausdorff space, $E$ a normed linear space,
	$M$ a subset of $\C(X,\R)$ which separates the points of $X$,
	$\vi\in E\setminus\{0\}$
	and $L_n$ a sequence of monotonically regular operators of $\C(X,E)$.
	If \[L_n({\bf h})\rightarrow {\bf h}\ \mathrm{uniformly\ on\ }X\]
	for any ${\bf h}\in \{\uni\vi\}\cup \{ h^j \vi |\ h\in M,\ j=1,2\}$,
	then \[L_n(F)\rightarrow F\ \mathrm{uniformly\ on\ }X\]
	for any $F\in\C(X,E)$.
\et

Since $L^1_n$ is monotonically regular and $C_\infty$ is compact, we shall use Theorem \ref{apt}
to prove our Theorem \ref{cube}.

\medskip

\bp{of Theorem \ref{cube}}
For $j\ge 1$, let $\pr_j:C_\infty\to \R$ be the canonical projection: $\pr_j(t)=t_j$.
Let $\vi\in E$ be a non zero constant, since
 $M=\{\pr_j | j\ge1\}$ separates the points of $C_\infty$, 
it is sufficient to check the
convergences on the test function: $\uni\vi$, $\pr_j\vi$ and $\pr_j^2\vi$. 
\bern L_n^1(\uni\vi)(t)&=&\sum_{h\in A_n}\vi \varphi_{n,h}(t)\\
		&=&\vi \sum_{h_1=0}^n{n \choose h_1}t_1^{h_1}(1-t_1)^{n-h_1} \cdots
		\sum_{h_n=0}^n{n\choose h_n} t_n^{h_n}(1-t_n)^{n-h_n}=\vi.
\eern
For $j>n$,
\bern L^1_n(\pr_j\vi)(t)&=&\sum_{h\in A_n}\pr_j
		\left(\frac h n\right)\vi \varphi_{n,k}(t)=0, \\
	L^1_n(\pr_j^2\vi)(t)&=&\sum_{h\in A_n}\pr_j^{2}
		\left(\frac h n\right)\vi\varphi_{n,k}(t)=0, \eern
while for $j\le n$,
\bern L^1_n(\pr_j\vi)(t)&=&\sum_{h\in A_n}\pr_j\left(\frac h n\right)\vi\varphi_{n,k}(t)=
		\vi \sum_{h\in A_n}\frac{h_j}{n}\varphi_{n,k}(t)\\
	&=&\vi \sum_{h_1=0}^n{n \choose h_1}t_1^{h_1}(1-t_1)^{n-h_1} \cdots
		\sum_{h_j=0}^n\frac{h_j}{n}{n\choose h_j} t_j^{h_j}(1-t_j)^{n-h_j}\\
	&&\qquad\qquad\cdots \sum_{h_n=0}^n{n\choose h_n} t_n^{h_n}(1-t_n)^{n-h_n}=t_j\vi,\\
	L^1_n(\pr_j^2\vi)(t)&=&\vi \sum_{h\in A_n}\frac{h_j^2}{n^2}\varphi_{n,k}(t)\\
	&=&\vi \sum_{h_1=0}^n{n \choose h_1}t_1^{h_1}(1-t_1)^{n-h_1} \cdots
		\sum_{h_j=0}^n\frac{h_j^2}{n^2}{n\choose h_j} t_j^{h_j}(1-t_j)^{n-h_j}\\
	&&\qquad\qquad\cdots \sum_{h_n=0}^n{n\choose h_n} t_n^{h_n}(1-t_n)^{n-h_n}=
			 t_j^2\vi+ \frac{t_j-t_j^2}{n}\vi.
\eern
From these identities, we conclude the proof.
\ep

\subsection{Hilbert case: proofs}

We begin recalling the definition (cfr. \cite{ld00, prolla})
\bd Let $L:D(L)\to\F(X,E)$, $S:D(S)\to\F(X,\R)$ be linear operators, with $D(L)$ and
	$D(S)$ subspaces of $\F(X,E)$ and $\F(X,\R)$,  respectively. $L$ is said to be
	\emph{dominated by} $S$ if 
$$\norm F\in D(S) \ \ {\mathrm{and} \ \ }	\norm{L(F)(t)}\le S(\norm{F})(t)\]
	for any $F\in D(L)$ and $t\in X$.
\ed

As already stated in Remark \ref{inter}, the operators acting
on vector-valued mappings and on scalar functions will be denoted
with the same symbol $L_n^2$.
Therefore, the operator $L_n^2:\F(\Gamma,E)\to \F(\Gamma,E)$ is
dominated by $L_n^2:\F(\Gamma,\R)\to \F(\Gamma,\R)$.

In order to prove Theorem \ref{hilbert} and \ref{ucb}, we shall use
the results stated in \cite{ld00}, which, for sake of completeness,
we report below.

\bt\label{teo:shimond.vett} Let $Y$ and $E$ be normed spaces,
	$X$ be a convex subset of $Y$, $K\subset X$  and for any $n\ge 1$,
	$L_n\colon D(L_n)\to\F(K,E)$ be a $S_n$-regular 
	linear operator dominated by the positive linear operator
	$S_n\colon D(S_n)\to\F(K,\R)$. We suppose that, for every $n\ge1$,
	$UCB(X,E)\subset D(L_n)$, $UCB(X,\R)\subset D(S_n)$
	and $\psi_t^2\decl\norm{\cdot -t}^2\in D(S_n)$
	for some (and hence for all) $t\in Y$.
	Then for each $F\in UCB(X,E)$, $t\in K$ and $\delta>0$, one has
	\be \norm{L_n(F)(t)-F(t)} \le
	\norm{F(t)} \abs{S_n(\uni)(t)-1}+\omega(F,\delta)\left[S_n(\uni)(t)+
		\delta^{-2}\gamma_n^2(t)\right],\!\label{dis:shimond.vett} \ee
	where $\gamma_n^2(t)\decl S_n(\psi_t^2)(t)$.
\et

From Theorem 4.1 and Remarks 4.2 and 4.3 in \cite{ld00}, we deduce the following

\bt\label{m} Let $Y$ be a real reflexive Banach space, $E$ normed space,
	$X$ a convex subset of $Y$ closed and unbounded or open,
	$K$ a bounded, closed convex subset
	of $X$ and $g:X\to\R$ satisfying the following conditions:
	$g$ is strictly positive, strictly convex, Fr\'echet differentiable
	on $K$, $g'(K)$ is bounded in $Y'$ and the function
	\[ h(t,u)\decl g(u)-\left[g(t)+\langle g'(t),u-t\rangle\right],\]
	is lower semicontinuous with respect to weak topology. 
	Moreover, setting $B_n\decl g^{-1}([0,n])$, we
	require that $K\subset B_n$, $B_n$ is bounded, $X\setminus B_n\neq \emptyset$ and
	\[ \lim_{\stackrel{\scriptstyle \norm t\rightarrow\infty}{t\in X}}
		\frac{g(t)}{\norm{t}}=+\infty.\label{eqn:ipo.cresc.inf.g}\]
	For each $n\ge1$, let $L_n\colon D(L_n)\to\F(K;E)$
	be a $S_n$-regular linear operator dominated by the linear positive operator	
	$S_n\colon D(S_n)\to\F(K,\R)$, with $\K(X,E,g)\subset D(L_n)$, $\K(X,\R,g)\subset D(S_n)$ and $g,h\in D(S_n)$.
	If for every continuous linear functional $\phi\in Y'$, the convergences
	\be S_n(\uni)(t)\rightarrow 1,\ 
		S_n(\phi_{|_X})(t)\rightarrow \phi(t)\ \mathit{and\ }S_n(g)(t)\rightarrow g(t)\label{test}\ee
	hold uniformly for $t\in K$, then for every $F\in\K(X;E,g)$ and $f\in\K(X,\R,g)$,
	\[L_n(F)(t)\rightarrow F(t)\ \mathit{and\ } S_n(f)(t)\rightarrow f(t)
		\mathrm{\ uniformly\ for\ }t\in K.\]
\et

In our case $Y$ is the real separable Hilbert space $H$, $X$ is the
set $\Gamma$ that results to be convex, unbounded and closed.
In order to prove the pointwise convergence in 
Theorems \ref{hilbert} and \ref{ucb}
we have only to check the convergences in (\ref{test}) and to evaluate
the quantities involved in (\ref{dis:shimond.vett}).
The proof of the uniform convergence will need of the following lemma.

\bl\label{lem1} Let $C\subset \ell^2$ be relatively compact. Then for any $\epsilon>0$,
	there exists an integer number $i=i(\epsilon,C)$, such that
	for every $x\in C$, we have
	\[ \sum_{j\ge i} x_j^2 <\epsilon.\]
\el
\bp{} Suppose, contrary to our claim, that there exist
	$\epsilon>0$ and a sequence $(x^i)_{i\ge 1}$ in $C$, such that
	\[ \sqrt{\sum_{j\ge i} (x^i_j)^2} \ge\sqrt{ \epsilon}\]
	for every $i\ge 1$.
	From the relatively compactness of $C$, there exists 
	$\bar x\in\overline C$ such that (up to a subsequence),
	$x^i\rightarrow\bar x$ (as $i\rightarrow \infty$).
	Thus, we have
\[ \sqrt\epsilon\le\sqrt{  \sum_{j\ge i} (x^i_j)^2} \le 
	\sqrt{\sum_{j\ge i} (x^i_j-\bar x_j)^2}+ \sqrt{  \sum_{j\ge i} (\bar x_j)^2}\le
        \norm{x^i-\bar x}+\sqrt{  \sum_{j\ge i} (\bar x_j)^2}, \]
	for every $i\ge 1$.
	Letting $i\rightarrow \infty$, we have a contradiction.
\ep


\bp{of Theorem \ref{hilbert}}
We begin fixing $A\subset \Gamma$ relatively compact and set $K$ the
compact convex hull of $A$.
Setting $g(u)\decl 1+\norm{u}^2$, 
we have that the function
\[h(t,u)=\norm{t}^2 +\norm{u}^2-2\langle t,u\rangle,\]
is lower semicontinuous for the weak topology.
Choosing $\lambda$ such that $K\!\subset\! g^{-1}([0,\lambda])$,
we have that the hypotheses of Theorem \ref{m} are satisfied.

Now, with the same computations of the proof of Theorem \ref{cube}, 
we evaluate the
convergences on the test functions.

We begin with
\be\label{eq:test0} L_n^2(\uni)(t)=\sum_{h\in A_n}\varphi_{n,h}(t)=1.\ee

Let us denote with
$(e_j)_{j\ge1}$ the dual base of $(a_j)_{j\ge 1}$ (that is the base of the dual space $H'$
such that $\langle e_i,a_j\rangle=\delta_{ij}$).
For $j>n$,
\bern L_n^2(e_j)(t)&=&\sum_{h\in A_n}e_j\left(\frac h n\right)\varphi_{n,k}(t)=0,\\
	L_n^2(e_j^2)(t)&=&\sum_{h\in A_n}e_j^{2}\left(\frac h n\right)\varphi_{n,k}(t)=0,\eern
while for $j\le n$,
\bern L_n^2(e_j)(t)&=&\sum_{h\in A_n}e_j\left(\frac h n\right)\varphi_{n,k}(t)=t_j,\\
L_n^2(e_j^2)(t)&=&\sum_{h\in A_n}\frac{h_j^2}{n^2}\varphi_{n,k}(t)=t_j^2+\frac{t_j-t_j^2}{n}.\eern

Let $t\in H$ and  $\phi\in H'$,  representing them as
$t=\sum_{j=1}^\infty t_ja_j$ and $\phi=\sum_{j=1}^\infty\phi_je_j$,
we have $\phi(t)=\sum_{j=1}^\infty\phi_jt_j$. Computing
\[ L_n^2(\phi)(t)=L_n^2(\sum_{j=1}^\infty\phi_je_j)(t)=\sum_{j=1}^\infty\phi_jL_n^2(e_j)(t)
	=\sum_{j=1}^n\phi_jt_j, \]
we obtain the convergence of $L_n^2(\phi)$ to $\phi$, uniformly on
bounded subsets of $\Gamma$.

Noting that $\psi_t^2(u)=\norm u^2 + \norm t^2 -2\langle t,u\rangle$,
in order to conclude the proofs, we have to evaluate
$L_n^2(\norm{\cdot}^2)$ on relatively compact subsets.
From identity
\[ \norm{t}^2=\sum_{j=1}^\infty t_j^2=\sum_{j=1}^\infty e_j^2(t), \]
we have
\[ L_n^2(\norm{\cdot}^2)(t)=\sum_{j=1}^\infty L_n^2(e_j^2)(t)
			=\sum_{j=1}^n\left( t_j^2+\frac{t_j-t_j^2}{n}\right), \]
and hence
\begin{eqnarray} L_n^2(\norm{\cdot}^2)(t)-\norm{t}^2&=&-\sum_{j>n} t_j^2
		+\sum_{j=1}^n\frac{t_j}{n}-\frac 1 n \sum_{j=1}^nt_j^2,\label{testg}\\ 
	L_n^2(\psi_t^2)(t)&=&\sum_{j>n} t_j^2
		+\sum_{j=1}^n\frac{t_j}{n}-\frac 1 n \sum_{j=1}^nt_j^2.\label{testpsi2}\end{eqnarray}

For the second term in the right hand side of (\ref{testg}) and 
(\ref{testpsi2}), the following estimate holds
\[ \sum_{j=1}^n\frac{t_j}{n}\le
		(\sum_{j=1}^n\frac{1}{n^2})^{1/2}  (\sum_{j=1}^n {t_j}^2)^{1/2}
\le \frac{1}{\sqrt n}\norm{t}.\]
Thus, the last two  terms in (\ref{testg}) and (\ref{testpsi2}) decay
 to 0 uniformly on bounded  subsets of $\Gamma$.
Therefore, the estimates and the convergences hold pointwise as claimed in
Theorem \ref{hilbert}. 
The uniform convergences on $K$ (and hence on $A$) follow from
the uniform convergence of $\sum_{j>n}t_j^2$ to 0, and this is stated
in the above Lemma (\ref{lem1}).
\ep

\bp{of Theorem \ref{ucb}} In order to prove the estimates in the statement of 
  Theorem \ref{ucb}, taking into account  Theorem \ref{teo:shimond.vett},
  it is sufficient to compute $L^2_n(\uni)(t)$ and 
  $L^2_n(\norm{\cdot -t}^2)(t)$.
  These quantities are already computed in the proof of Theorem~\ref{hilbert}.
  Hence from (\ref{eq:test0}) and (\ref{testpsi2}), we obtain the stated estimates.
\ep


\medskip
\bp{of Proposition \ref{lip}}
The preserving properties of Proposition \ref{lip}, follow
from the definition of $L_n^2$ and from Theorem \ref{tb}.
For instance, the inclusion
$L_n^2(Lip_M(\Gamma))\subset Lip_{\sqrt nM}(\Gamma)$
follows from 2. of Theorem \ref{tb} and the relation
\be \sum_{i=1}^n \abs{t_i}\le \sqrt n\left(\sum_{i=1}^n \abs{t_i}^2\right)
	\le \sqrt n \norm{t}. \label{norm}\ee

See Proposition \ref{pro2} for more general cases.\ep

When $E$ is the real line, it remains to prove that the analogues of
properties 4. and 5. of Theorem \ref{tb} fail for $L_n^2$.
Indeed, it is enough to consider what happens with the functionals $e_j$,
the base of $H'$: for $j>n$, $L_n^2(e_j)=0$ and
$n\ge j$, $L_n^2(e_j)=e_j$.
Thus, one can conjecture that the properties hold definitively, that is,
for any $f\in \K(\Gamma,\R,g)$ convex, there exists an integer $\nu$
such that, for $n\ge\nu$,
$L_n^2(f)\ge f$ and $L_n^2(f)\ge L_{n+1}^2(f)$.
Tough, even this conjecture is doomed to fail.
Indeed, let $\bar f$ be the function defined as 
$\bar f\decl \sum_{j\ge 1}\frac{e_j}{2^j}$.
The function $\bar f$ is convex and belongs to $\K(\Gamma,\R,g)$.
Computing
\[L_n^2(\bar f)(t)=\sum_{j=1}^n\frac{1}{2^j}\left(t^2_j+\frac{t_j-t_j^2}{n}\right),\]
and applying at $\bar t=(1,\dots,1,t_{n+1},t_{n+2},\dots)$, we obtain
\bern \bar f(\bar t)-L_n^2(\bar f)(\bar t)&=&\sum_{j\ge n+1}\frac{t^2_j}{2^j}\ge 0,\\
L_{n+1}^2(\bar f)(\bar t)-L_n^2(\bar f)(\bar t)&=&
\frac{1}{2^{n+1}}\left(t^2_{n+1}+\frac{t_{n+1}-t_{n+1}^2}{n+1}\right)\ge0,\eern
that prove our claims.

\section{A generalization}
In this section we generalize the proposed scheme.
We start with a generic sequence of positive linear operators,
and as in Bernstein polynomials case,
we obtain approximation processes for vector-valued mappings
defined on  subsets of an infinite dimensional Hilbert  space.

Let $E$ be a Banach space, $I$ a Hausdorff space, $J\subset I$ and
for $n\ge 1$, and $t\in J$, $\mu_n(\cdot;t)$  a probability measure
on $\sigma$-algebra of all Borel subset of $I$. 
With $L^1(I,E,\mu_n(\cdot;t))$, we denote the subspace of $\F(I,E)$ of
all $\mu_n(\cdot;t)$-integrable functions.
We consider the linear integral operator
$L_{n,1}:L^1(I,E,\mu_n(\cdot;\cdot))\to \F(J,E)$, defined as
\[L_{n,1}(f)(t)\decl \int_I f(u) d\mu_n(u;t).\]
From the measure $\mu_n(\cdot,t)$, we define for $n\ge 1$, $k\ge 1$
and $t=(t_1,\cdots,t_k)\in J^k$ the product measure
$\mu_{n,k}(\cdot;t)\decl\bigotimes_{i=1}^k\mu_n(\cdot;t_i)$, and then
we consider the associated integral operator:
\[L_{n,k}(f)(t)\decl\int_{I^k}f(u)d\mu_{n,k}(u;t)= \int_{I^k} f(u_1,\dots,u_k) d\mu_n(u_1;t_1)\otimes \cdots
	\otimes d\mu_n(u_k;t_k),\]
for $t=(t_1,\dots,t_k)\in J^k$, and 
$f\in L^1(I^k,E, \mu_{n,k}(\cdot;t))$

We fix $s=(s_i)_{i\ge1}\in I^{\nat^*}$.
For $f:I^{\nat^*}\to E$, the symbol $f_k$ stands for the function
$f_k: I^k\to E$ defined as
$f_k(t_1,\dots,t_k)\decl f(t_1,\dots,t_k,s_{k+1},s_{k+2},\dots)$.
In the other direction, for $f:I^k\to E$, the symbol $\tilde f$ denotes
the function $\tilde f:I^{\nat^*}\to E $, defined as
$\tilde f(t)\decl f(t_1,\dots,t_k)$.
Finally, for  $f:I^{\nat^*}\to E$ such that
$f_n\in L^1(I^n,E, d\mu_{n,n}(\cdot;t))$, for any $t\in J^{\nat^*}$,
we define 
\[L_n(f)\decl (L_{n,n}(f_n))\tilde{}.\]
It is immediate to check that $L_n$ is a monotonically regular operator.

One can hope that some property of $L_{n,1}$ are inherited
from $L_n$. For instance, choosing $L_{n,1}=B_{n,1}$,
the Bernstein operators, $s=0$, it results
$L_n^2(f)=(B_{n,n}(f_n))\tilde{} $, for $f:\Gamma\to E$.
If we  define $L_n^2$ with a generic $s\in\Gamma$,
\[ L^2_n(f)(t)\decl \sum_{h\in A_n} f\left(\frac{h_1}{n}, \dots,\frac{h_n}{n},
	s_{n+1},s_{n+2},\dots \right)\varphi_{n,h}(t), \]
then this variation is not essential. Indeed,
Theorems \ref{hilbert}, \ref{ucb} and their proofs are the same, and 
with a small change of the function $\bar f$, one can
show that analogue properties of 4. and 5. of Theorem \ref{tb}
do not hold.

\bt\label{gat} In the same setting of subsection 2.2 and with
	the above notation, let $I=J$ be a real interval with $0\in I$,
	$\Gamma'\decl\{t\in H|\ t_i\in I \}$, and fix
	$s=(s_i)_{i\ge 1}\in\Gamma'$.
	We assume that $e_2\in L^1(I,\R,\mu_n(\cdot;t))$ 
	for every $n\ge 1$ and $t\in J$,
	$L_{n,1}(e_1)=e_1$ and
	$L_{n,1}(e_2)=e_2+e_2o(1)+e_1o(\frac{1}{\sqrt n})+o(\frac{1}{n})$.
\begin{enumerate}
\item If $F\in\K(\Gamma',E, 1+\norm{\cdot}^2)$, or $F\in UCB(\Gamma',E)$, then
	\[L_n(F)\rightarrow F\]
	uniformly on relatively compact subsets of $\Gamma'$.
\item If $L_{n,1}(Lip_1(I))\subset Lip_1(I)$,
	then $L_{n}(Lip_M(\Gamma'))\subset Lip_{\sqrt n M}(\Gamma')$.
\end{enumerate}
\et

We note that the conditions of Theorem \ref{gat} are satisfied by
many operators, e.g. Sz\'asz--Mirakjan operators, Baskakov operators,
Post--Widder operators.

\boss	In the assumption
	$L_{n,1}(e_2)=e_2+e_2o(1)+e_1o(\frac{1}{\sqrt n})+o(\frac{1}{n})$,
	the last term cannot be substituted with the weaker condition $O(\frac 1n)$.
	Indeed, let $L_{n,1}$ be the Gauss-Weierstrass
	operators, defined for $t\in \R$ and $f\in \C(\R,\exp(e_2))$, as
\[ L_{n,1}(f)(t)\decl \sqrt{\frac{n}{\pi}}\int_{\RR} f(u) e^{-n(u-t)^2} du.\]
	It results $L_{n,1}(\uni)=\uni$, $L_{n,1}(e_1)=e_1$,
	$L_{n,1}(e_2)=e_2+\frac{1}{2n}$ and
	$L_{n,1}$ approximates uniformly on bounded sets the functions
	belonging to $\C(\R,\exp(e_2))$ (see \cite{shaw80}).

	Choosing $s_i=0$, with same notation as before,
	we get $L_{n}(\norm{\cdot}^2)(t)=\sum_{i=1}^n t_i^2 + 1/2$,
	which converges to $\norm{t}^2 +1/2$. 
	Therefore, we cannot conclude that $L_n$ is an approximation process for functions
	belonging to $\K(H,\R,1+\norm{\cdot}^2)$.
\eoss

In order to prove the preserving property 2., we give the following result.

\bpr\label{pro2} Let $(I,d)$ be metric space. Consider the metric space $I^k$
	endowed with the distance
	$d_k(t,\tau)\decl\sum_{i=1}^k d(t_i,\tau_i)$.
	If $L_{n,1}(Lip_1(I))\subset Lip_1(I)$, then for any $k\ge 1$
	$L_{n,k}(Lip_1(I^k))\subset Lip_1(I^k)$.
\epr
\bp{} We  shall argue  by induction on $k$.
	For $k=1$, the property holds by hypothesis.
	We assume that it is true for $k-1$.
	Let $f\in Lip_1(I^{k})$, $t,\tau\in I^{k-1}$ and $t_{k},\tau_k\in I$.
	Using the integral nature of the operators $L_{n,k}$, one gets
\bern \lefteqn{L_{n,k}(f)(t,t_k)-L_{n,k}(f)(\tau,\tau_k)=}\\
	&&=\int_I\big[ L_{n,k-1}(f(\cdot,u_k))(t)-
		L_{n,k-1}(f(\cdot,u_k))(\tau)\big]d\mu_n(u_k;t_k)\\
	&&\qquad+\int_{I^{k-1}}\left[ L_{n,1}(f(u,\cdot))(t_k)- L_{n,1}(f(u,\cdot))(\tau_k)\right]d\mu_{n,k-1}(u;t).\eern
Thus, since $f_{|I^{k-1}}\in Lip_1(I^{k-1})$, we obtain
\bern \lefteqn{\norm{L_{n,k}(f)(t,t_k)-L_{n,k}(f)(\tau,\tau_k)}}\\
	&&\le \int_I d_{k-1}(t,\tau) d\mu_n(u_k;t_k)
	+\int_{I^{k-1}} d(t_k,\tau_k) d\mu_{n,k-1}(u;t) \\
	&&\qquad\qquad\qquad\qquad\qquad\qquad=d_k((t,t_k),(\tau,\tau_k)),\eern
which allows us to conclude the proof of the proposition.
\ep

\bp{of Theorem \ref{gat}}
The proof of the approximation property of Theorem \ref{gat},
using Theorems \ref{teo:shimond.vett} and \ref{m}, is the same of
Theorems  \ref{hilbert} and \ref{ucb}.

In the setting of Theorem \ref{gat},
the inclusion $L_n(Lip_M(\Gamma'))\subset Lip_{\sqrt nM}(\Gamma')$
is now immediate.
Indeed, if $f\in Lip_1(\Gamma')$, then also its restriction
$f_{k}$ belongs to $Lip_1(I^k)$, for every $k\ge 1$.
Hence, $L_{n,n}(f_n)\in Lip_1(I^n)$, and from inequality (\ref{norm}),
we get the thesis.
\ep







\end{document}